\documentclass[%
 reprint,
 amsmath,amssymb,
 aps,
]{revtex4-1}

\usepackage{color}
\usepackage{graphicx}% Include figure files
\usepackage{dcolumn}% Align table columns on decimal point
\usepackage{bm}% bold math
\usepackage{hyperref}% add hypertext capabilities
\begin{document}

\preprint{APS/123-QED}

\title[]{Shape Coherence and Finite-Time Curvature Evolution}%:\\with Forced Linebreak\footnote{Error!}}% Force line breaks with \\
%\thanks{Footnote to title of article.}

\author{Tian Ma}
\author{Erik M. Bollt}%
 \email{bolltem@clarkson.edu}
\affiliation{ 
Department of Mathematics and Computer Science, Clarkson University, USA.%\\This line break forced with \textbackslash\textbackslash
}%

\date{\today}% It is always \today, today,

\begin{abstract}
We introduce a definition of finite-time curvature evolution along with our recent study on shape coherence in nonautonomous dynamical systems. 
Comparing to slow evolving curvature preserving the shape, large curvature growth points reveal the dramatic change on shape such as the folding behaviors in a system. 
Closed trough curves of low finite-time curvature (FTC) evolution field indicate the existence of shape coherent sets, and troughs in the field indicate most significant shape coherence. Here we will demonstrate these properties of the FTC, as well as contrast to the popular Finite-Time Lyapunov Exponent (FTLE) computation, often used to indicated hyperbolic material curves as Lagrangian Coherent Structures (LCS).  We show that often the FTC troughs are in close proximity to the FTLE ridges, but in other scenarios the FTC indicates entirely different regions.
%\begin{description}
%\end{description}
\end{abstract}

\pacs{Valid PACS appear here}% PACS, the Physics and Astronomy
\maketitle

%\tableofcontents

%%%%%%%%%%%%%%%%%%%%%%%%%%%%%%%%%%%%%%%%%%%%%%%%%%%%%%%%%%%%%%%%%%%%%%%%%%%%%%%%%%%%%%%%%%%%%%%
%       Intro
%%%%%%%%%%%%%%%%%%%%%%%%%%%%%%%%%%%%%%%%%%%%%%%%%%%%%%%%%%%%%%%%%%%%%%%%%%%%%%%%%%%%%%%%%%%%%%%
%{\{red}Intro}
Coherence has clearly become a central concept of interest in nonautonomous dynamical systems, particularly in the study of  turbulent flows, 
with many recent papers designed toward describing,  quantifying and constructing such sets. \cite{MB1, HB, FSM, MB, FK, OG2, TR}.  There have been a wide range of notions of coherence, from spectral, \cite{LumleyHolmes}, to set oriented, \cite{DJ} and through transfer operators \cite{FK,FSM} as well as variational principles \cite{M}, and even topological methods, \cite{Thiffeault2, Ross}.  
% A general perspective of set oriented analysis of coherence seems to emphasize a discussion of transport.  
 Traditionally there has been an emphasis on vorticity \cite{HU}, but generally an understanding that, coherent motions have a role in maintenance (production and dissipation) of turbulence in a boundary layer, \cite{RB}.
A number of theories have been developed to model and analyze the dynamics in the Lagrangian perspective (moving frame), 
such as the geodesic transport barriers \cite{HB} and transfer operators method \cite{FSM}.  
These have included analysis of coherence in important problems such as how regions of fluids are isolated from each other \cite{Thiffeault2} 
including in prediction of oceanic structures \cite{FK1} and atmospheric forecasting 
\cite{Ross1, Ross2}, 
especially for the understanding of movement of pollution including such as oil spills, 
\cite{Haller3, Mezic, BolltGulf}.
Whatever the perspectives taken, we generally interpretively  summarize that coherent structures can be taken as a region of simplicity, within the observed time scale and  stated spatial scale, perhaps embedded within an otherwise possibly turbulent flow, \cite{HB, FSM, FK, MB1}.  
%More about how coherence is timely and important - in intro and conclusion.

 In particular, the ridges from Finite Time Lyapunov Exponents (FTLE) fields have been widely used  \cite{H1, H2,SLM,TR} to indicate hyperbolic material curves, often called Lagranian coherent structures (LCS). We contrast here the fundamental nonlinear notions of ``stretching" encapsulated in the FTLE concept to ``folding" which is a  complementary concepts of a nonlinear dynamical system which must be present if a material curve can stretch indefinitely within a compact domain.  As we will show that exploring the much-overlooked folding concepts leads to developing curvature changes of material curves yielding an elegant description of coherence that we call shape coherence \cite{MB}. We introduce here a method of visualizing propensity  of a material curve to change its curvature, which we call the Finite-Time Curvature (FTC) field.  
Contrasting the FTC to the FTLE, we will illustrate that  sometimes the FTC troughs indicative of shape coherence are often co-mingled in close proximity to ridges of the FTLE, and in such case they indicate a generally similar story.  However we show that in many cases the FTC troughs occur in locations not near an FTLE ridge, indicating entirely different regions.  Thus we view these as complementary concepts, stretch and fold, as revealed by the traditional FTLE and the here  introduced FTC.
%On the other hand, there are several studies of relationships between folding behaviors and curvature, \cite{Thiffeault, Thiffeault1,TA}
%See discussion of curvature in the context of dynamical systems also in \cite{Thiffeault, Thiffeault1, TA, OG, OG1, OG2, OG3, HB}

We have recently presented a mathematical interpretation of coherence \cite{MB} in terms of a definition  of {\it shape coherent sets},
 motivated by a simple observation regarding sets that ``hold together'' over finite-time in nonautonomous dynamical systems. 
As a general setup, assume an area preserving system that can be represented, 
%\begin{equation}\label{jj}
$\dot{z}=g(z,t), $
%\end{equation}
for $z(t)\in \Re^2$, with enough regularity of  $g$ so that that a corresponding flow, $\Phi_T(z_0):\Re\times \Re^2\rightarrow \Re^2$, exists.  To capture the idea of a set that roughly preserves its own shape, we define \cite{MB} 
the {\it shape coherence factor} $\alpha$ between two sets $A$ and $B$ 
under an area preserving  flow $\Phi_t$ over a finite time interval $[0, T]$,
{\begin{eqnarray}\label{shapecoherenced}
\displaystyle \alpha(A, B,T) := \sup_{S(B)} \frac{m( S(B) \cap \Phi_{ T}(A)) }{m(B)},
\end{eqnarray}}
where $m(\cdot)$ here denotes Lebesgue measure, and we restrict the domain of $\alpha$ to sets such that $m(B)\neq 0$ by assumptions to follow that $B$ should be a fundamental domain \cite{Complex}. 
Here,
$S(B)$ is the group of transformations of rigid body motions of $B$, specifically translations and rotations descriptive of {\it frame invariance,}\cite{Carmo}.  
We say $A$ is finite time shape coherent to $B$ with shape coherence factor $\alpha$, under the flow $\Phi_T$ after the time epoch $T$.  We  call $B$ the {\it reference set}, and $A$ shall be called the {\it dynamic set}. 
If we choose $B = A$, we can verify to what degree a set $A$ preserves its shape over the time epoch $T$.
Notice that the shape of $A$ may vary during the time interval, 
but for a high shape coherence, the shapes must be similar at the terminal times. 
By the area preserving assumption, $0 \leq \alpha \leq 1$, and values closer to $1$ indicate a set for which the otherwise nonlinear flow restricted to $A$ is much simpler, at least  on the time scale $T$ and on the spatial scale corresponding to $A$; that is $\Phi_T|_A$, the flow restricted to $A$ is roughly much simpler than a turbulent system, as it is much more like a rigid body motion.  %Note that 
This does not preclude  on finer scales, that there may be  turbulence within a shape coherent set.

Recall that for any material curve, $\gamma(s,t)=(x_1(s,t),x_2(s,t))$ of initial conditions defining an initial segment $\gamma(s,0)=(x_1(s,0),x_2(s,0))$, $a\leq s\leq b$   where each point on the curve evolves in time $t$ according to the differential equation, the curvature at time $t$ may be written in terms of the parametric derivative along the curve segment, $d/ds:='$,
$k(s,t)=\frac{|x_1'x_2''-x_2'x_1''|}{(x_1'^2+x_2'^2)^{3/2}}$.
We will relate the pointwise changes of this curvature function for points on those material curves that correspond to shape coherence.  
%Note that this perspective of coherence is inherently different from other studies of coherence, most notably those based on finite time Lyapunov exponents (FTLE) \cite{H1, H2,SLM,TR} since both first and second derivative information are explicitly included here.

\begin{figure}[h]
  \centering
  \includegraphics[width=0.4\textwidth]{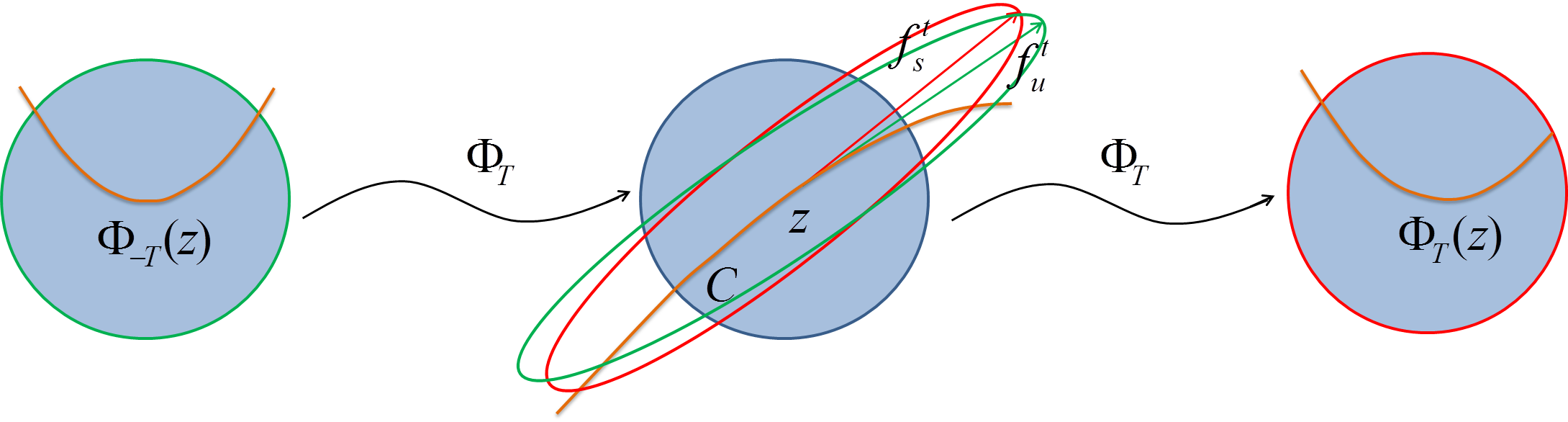} 
  
%  \includegraphics[width=0.4\textwidth]{Tang}
%  (b)
% % PLUS maybe a more true manifolds story and from computation too.
  \caption{
Tangency points of stable and unstable foliations are highlighted as at such points, infinitesimal curve elements experience curvature that  terminally} evolves slowly in the  time epoch. Notice that the stable foliation $f_s^T(z)$ is the major axis of preimage of variations from $\Phi_T(z)$ and correspondingly $f_u^T(z)$ is the major axis of image of variations from $\Phi_{-T}(z)$. 
  \label{IdeaOfFTC}
\end{figure}

\begin{figure}[h]
  \centering
  \includegraphics[width=0.48\textwidth]{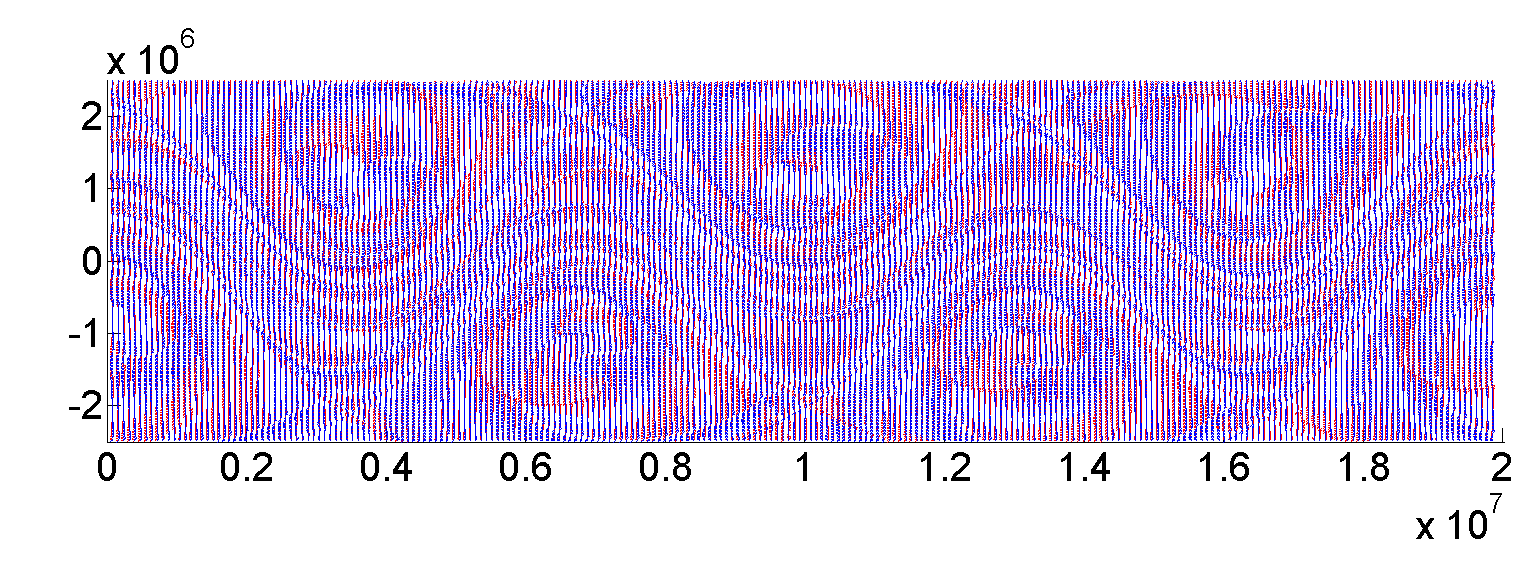} 
  (a) Foliation field.
  \includegraphics[width=0.46\textwidth]{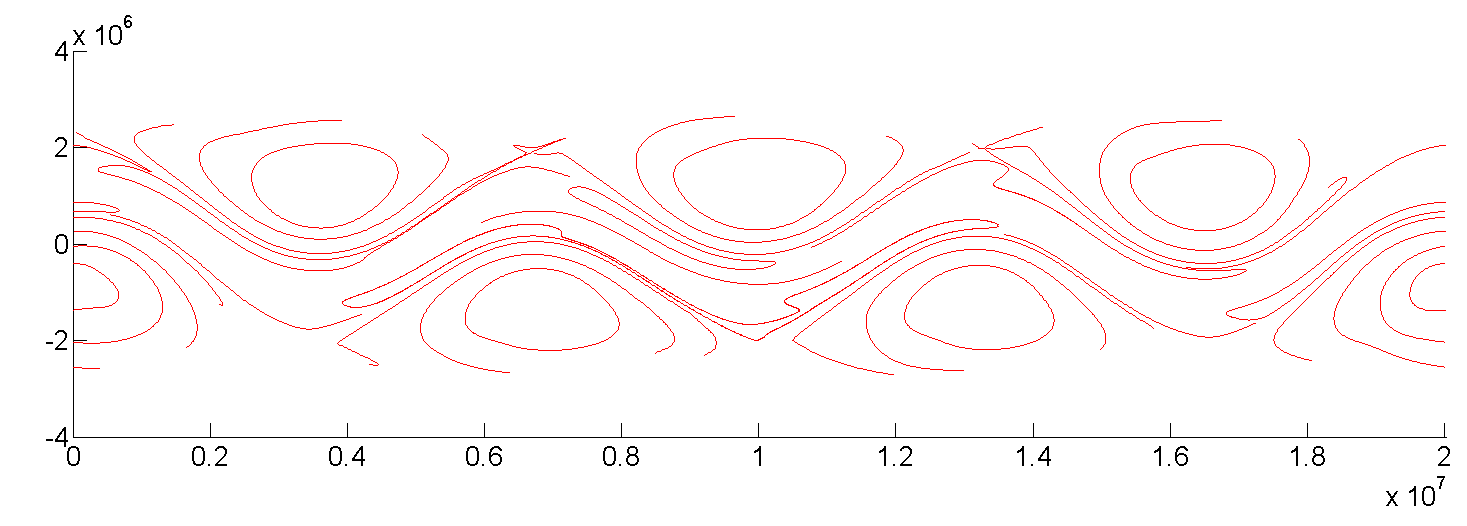} 
  (b) Nonhyperbolic splitting curves.
  \caption{(a) The finite-time stable and unstable foliation field $f_s^t(z)$ and $f_u^t(z)$ for the Rossby wave system. Notice that for each point, there are two vectors, and therefore an associated angle $\theta(z, t)$. (b)  zero-splitting curves corresponding to boundaries of shape coherent sets corresponding to a significant shape coherence factor, Eq.~(\ref{shapecoherenced}).
  The parameters are set, $U_0=44.31, c_2=0.2055U_0, c_3=0.462U_0, A_3=0.3, A_2=0.12, A_1=0.075$  \cite{RBB},
and time epoch $T=10 $ days. 
}
 \label{FoliationFieldAndCurves}
\end{figure}

 The analysis of the geometry of shape coherent sets $A$ 
  depends on the boundary of these sets, $\partial A$, which we restrict in  the following to  simply connected  sets such that the boundary  is a smooth and simple closed curve,  $\partial A=\gamma(s), 0 \leq s \leq 1$, and these are often called ``fundamental domains" \cite{Complex}.  These $B=A$ are in the domain of $\alpha$.  We may relate shape coherence to the classical differential geometry  whereby two curves are defined to be congruent if their underlying curvature functions can be exactly matched, pointwise, \cite{Carmo}.  Therefore, considering the Frenet-Serret formula \cite{Carmo},  it can be  proved \cite{MB1} through a series of regularity theorems that those sets with a slowly  evolving propensity to change curvature correspond to boundaries of sets with a significant degree of shape coherence.  That is $\alpha(A,\Phi_T(A))\approx 1$.
Furthermore, a  sufficient condition theorem connects  geometry that  points $z$ where there is a tangency between finite time stable and unstable foliations $f_u^t(z)$, $f_s^t(z)$ must correspond to slowly changing curvature.  In Fig.~\ref{IdeaOfFTC}, we indicate the   geometry of stable and unstable foliations that correspond to tangency or near tangency where curves passing through such points experience slowly changing curvature, and hence indicative of points on the boundaries of shape coherent sets \cite{MB1}.    Hence, to find shape coherent sets lead us to the search for curves of tangency points as the boundaries of such set which we review below. Much has been written about   the role of how  stable and unstable manifolds can become reversed at tangency points in that errors can grow transversely to the the unstable manifolds as noted in  \cite{Kantz, BSLZ, ZB}.
Scaling relationships for frequency of given curvatures in \cite{Thiffeault, Thiffeault1, TA},  \cite{Liu,Drummond1,Drummond2,Pope,Ishihara}, as well as the propensity of curvature growth in turbulent systems \cite{OG,OG1,OG2,OG3} have both been studied.

We review, the finite time stable foliation $f_s^t(z)$ at a point $z$ describes the dominant direction of local contraction in forward time, and the finite time unstable foliation $f_u^t(z)$ describes the dominant direction of contraction in ``backward" time, and these vectors have a long history in the stability analysis of a dynamical system, particularly related to Lyapunov exponents and directions, \cite{Ulrich, BN}, and lately in \cite{HB}. See Fig.~\ref{IdeaOfFTC}.
The derivative,
$D\Phi_t(z)$ of the flow $\Phi_t(\cdot)$ evaluated at the point $z$ maps a circle  onto an ellipse, as does any general matrix, \cite{Geometry}
the infinitesimal geometry of a small disc of variations from near $\Phi_t(z)$  shown in Fig.~\ref{IdeaOfFTC}.  Likewise, a disc centered on $\Phi_t(z)$ pulls back under $D\Phi_{-t}({\Phi_{t}(z)})$ to an ellipsoid centered on $z$.  The major axis of that infinitesimal ellipsoid defines  $f_s^t(z)$, the stable foliation at $z$.  Likewise, from $\Phi_{-t}(z)$, a small disc of variations pushes forward under $D\Phi_{t}({\Phi_{-t}(z)})$ to an ellipsoid, the major axis of which defines, $f_u^t(z)$.
These major axis can be readily computed in terms of the  singular value decomposition  \cite{VanLoan} of derivative matrices, as noted regarding the Lyapunov directions \cite{Ulrich, BN, Oseledets} and recently \cite{DanielarXiv}.  Let,
$D\Phi_{t}({z})=U\Sigma V^*$, 
where $^*$ denotes the transpose of a matrix.  $U$ and $V$ are orthogonal matrices, and $\Sigma=diag(\sigma_1,\sigma_2)$ is a diagonal matrix.  
Indexing, $V=[v_1,v_2]$, and $U=[u_1,u_2]$, note that $D \Phi_{t}({z}) v_1=\sigma_1 u_1$ describes the vector $v_1$ at $z$ that maps onto the major axis, $\sigma u_1$ at $\Phi_t(z)$.  Since $\Phi_{-t}\circ \Phi_t(z)=z$, and $D\Phi_{-t}(\Phi_t(z)) D \Phi_t(z)=I$, then recalling  orthogonality of $U$ and $V$, yields,
$D\Phi_{-t}({\Phi_t(z)})=V\Sigma^{-1} U^*$, 
and $\Sigma^{-1}=diag(\frac{1}{\sigma_1},\frac{1}{\sigma_2})$.  Therefore, $\frac{1}{\sigma_2}\geq \frac{1}{\sigma_1}$, and the dominant axis of the image of an infinitesimal  disc from  $\Phi_t(z)$ comes from, 
$D \Phi_{t}({z}) u_2=\frac{1}{\sigma_2} v_2$.  
Hence,
\begin{equation}
 f_s^t(z)=v_2, \mbox{ and, } f_u^t(z)=\overline{u}_1,
 \end{equation}
where $v_2$ is the second right singular vector of $D\Phi_{t}({z})=U\Sigma V^*$ and likewise, $\overline{u}_1$ is the first left singular vector of 
$D\Phi_{t}({\Phi_{-t}(z)})=\overline{U} \mbox{ } \overline{\Sigma} \mbox{ } \overline{V}^*.$

% \cite{Thiffeault, Thiffeault1, TA, OG, OG1, OG2, OG3, HB}
 
%%%Following the above, for a given flow and time epoch, $\Phi_t$, the plane is doubly covered by two vector fields, $f_s^t(z)$ and $f_u^t(z)$.  Since one describes variations forward in time and one in backward in time, there is no reason there may not be tangencies.  
%%%When closed curves of points with such a tangency  can be found, then they must enclose significantly shape coherent sets.  Construction of these curves follows the implicit function theorem as applied to an angle function, 
%%%$ \theta(x, t): \Omega \times \mathbb{R}^+ \to [-\pi/2, \pi /2],$
%%%{\begin{eqnarray}
%%%\theta(z, t):=\arccos\frac{\left\langle f_s^t(z), f_u^t(z) \right\rangle}{\|f_s^t(z) \| \| f_u^t(z) \|},
%%%\label{AngleFunction}
%%%\end{eqnarray}}
%%%Thus we presented a constructive sufficient condition for significantly shape coherent sets \cite{MB1}.

%%%%%%%%%%%%%%%%%%%%%%%%%%%%%%%%%%%%%%%%%%%%%%%%%%%%%%%%%%%%%%%%%%%%%%%%%%%%%%%%%%%%%%%%%%%%%%%
%      FTC field
%%%%%%%%%%%%%%%%%%%%%%%%%%%%%%%%%%%%%%%%%%%%%%%%%%%%%%%%%%%%%%%%%%%%%%%%%%%%%%%%%%%%%%%%%%%%%%%
%{\color{red}Finite-Time Curvature}

\begin{figure}[htbp]
  \centering
  \includegraphics[width=0.98\linewidth, height = 0.1\textheight]{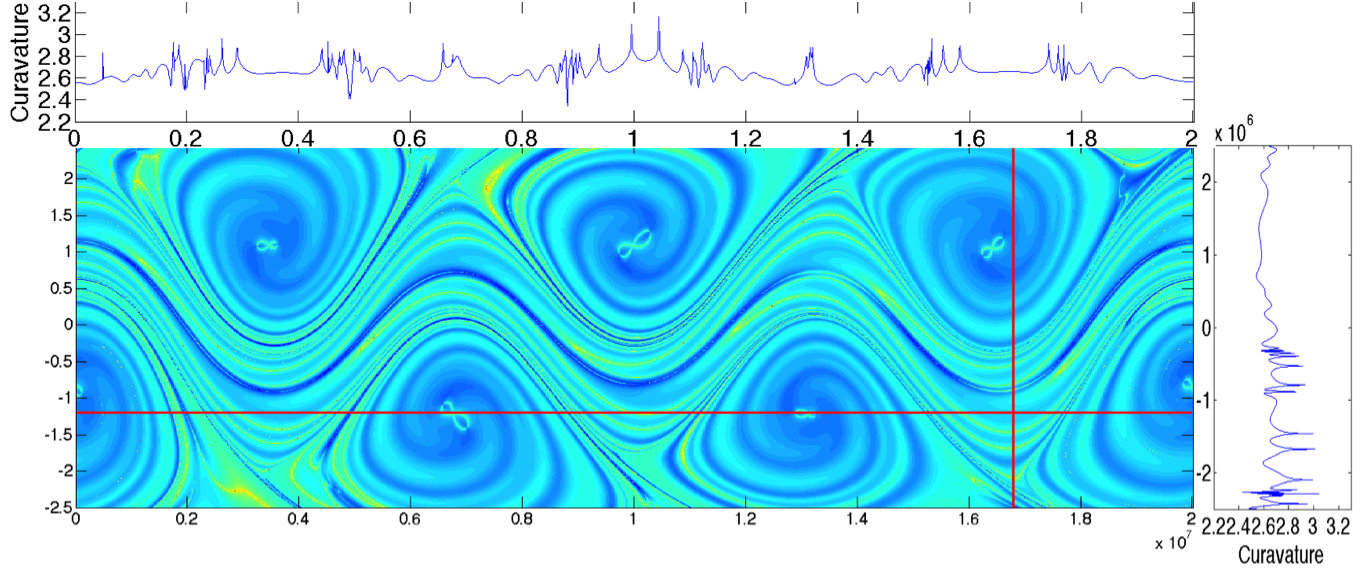} \\ % picture filename
  (a) maxFTC $U_0=44.31,c_3=0.461U_0, A_3=0.3$,
  \includegraphics[width=0.98\linewidth, height = 0.1\textheight]{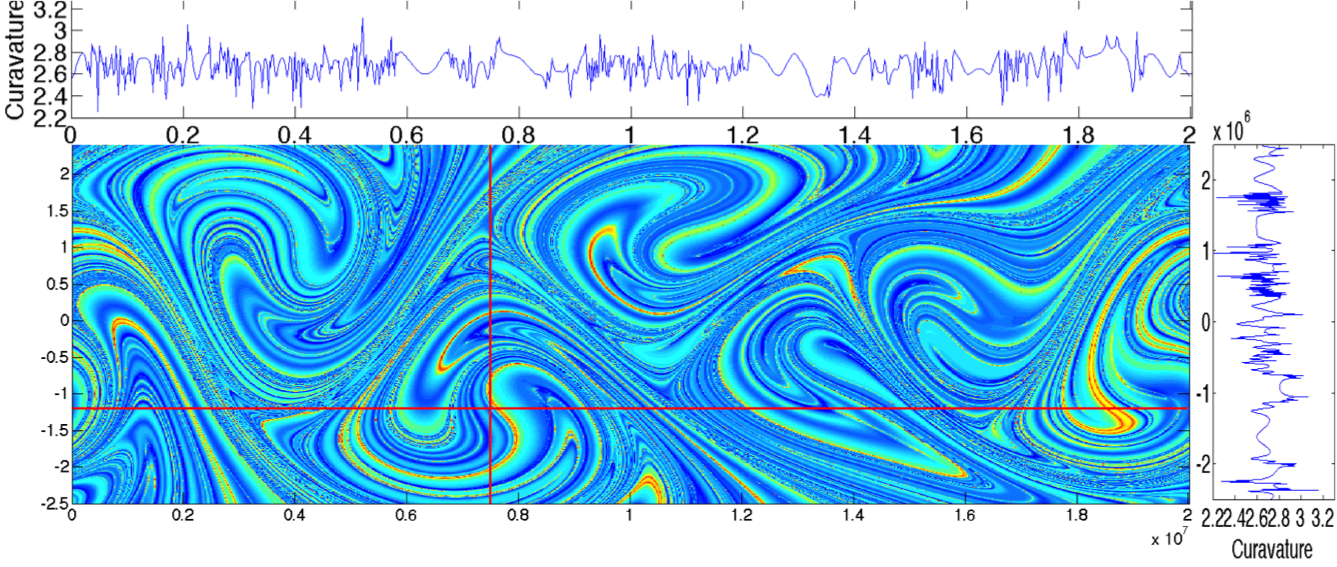}
  (b) maxFTC  $U_0=63.66, c_3=0.7U_0, A_3=0.2$.
%  \includegraphics[width=0.5\textwidth]{RWForPRLMin} \\ % picture filename
  %(c) minFTC $U_0=44.31,c_3=0.461U_0, A_3=0.3$,
   \includegraphics[width=0.5\textwidth]{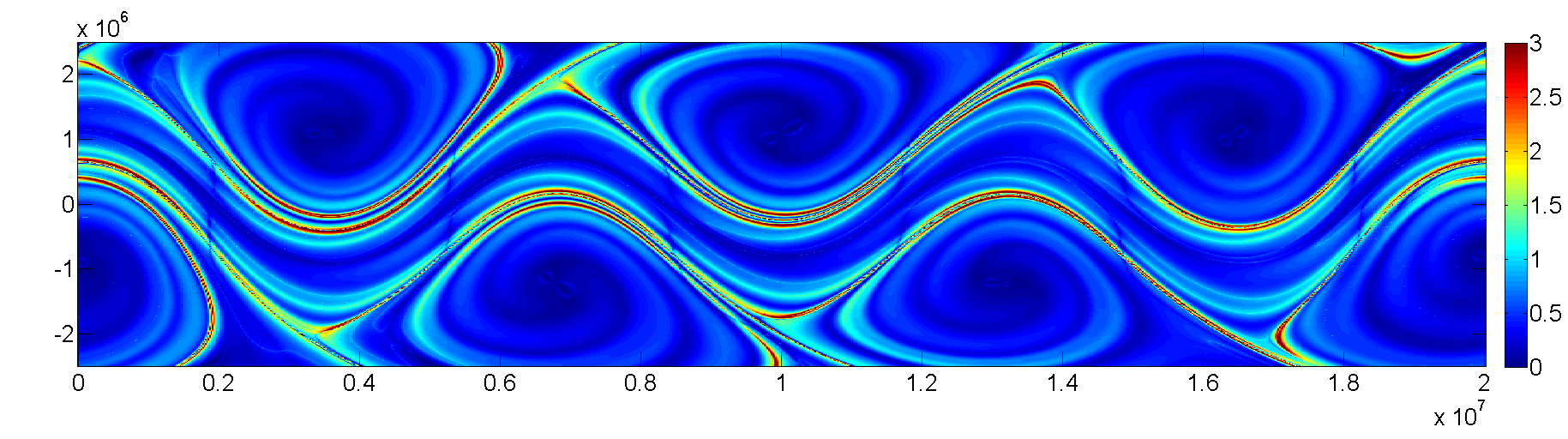}
  (c) FTC $U_0=44.31,c_3=0.461U_0, A_3=0.3$,
  \caption{(a) %and (b) show 
  shows the maxFTC fields $C_{t_0}^{t_0+\tau}(z)$ of different groups of parameters of the Rossby wave system. Note that in both figures, the level curves of relatively smaller maxFTC, $C_{t_0}^{t_0+\tau}(z)$ from Eq.~(\ref{ftc}), indicate that there exist material curves whose curvature changes slowly (blue curves) and these correspond to the zero-splitting curves in Fig.~\ref{FoliationFieldAndCurves}, \cite{MB1}. 
In the top and sides of (a) and (b)  we show a slice of the maxFTC function along the red lines shown respectively.  Large variation of these slice functions indicate the boundary of shape coherent sets; the interiors correspond to the slow variations of the function  and shape coherence, and generally those boundaries are indicated by low values of the maxFTC, small propensity to grow curvature.  The fast varying nature at boundaries indicates that high curvature change is often closely proximal to low curvature change, as indicated within (hetero)homoclinic tangle where tangencies and hyperbolicity often co-exist.
%(c) The minFTC field,  $c_{t_0}^{t_0+\tau}(z)$, Eq.~(\ref{ftcfield})   
(c) The FTC field,  $r_{t_0}^{t_0+\tau}(z)$, Eq.~(\ref{ftcfield}).
 % And large changes on curvature (red part) occur at the tips of the blue lines and the most twisted part of the system.
  }
  \label{RWFTC}
\end{figure}

%***probably move this to near the front***
%
%Here we  introduce a  simplified analysis and construction of shape coherence by  direct measurement and display of fields of  maximal rate of change of any possible curvature, which we  call the Finite Time Curvature field, (FTC).  The FTC allows us  to interpret sets of significant shape coherence, by direct  inspection of those points and curves corresponding to slowly evolving curvature,  with the interpretation cited above   theorems  connect shape coherence to slowly evolution of boundary curvature.
%
%****

For sake of further presentation, a specific example will be helpful. We choose  the Rossby wave \cite{RBB} system, an idealized zonal stratospheric flow.
Consider the Hamiltonian system
$dx/dt=-\partial H/ \partial y$,
$dy/dt=\partial H/ \partial x$, 
where
%\begin{eqnarray}
%\Phi(x,y,t)=c_3y&-&U_0Ltanh(y/L)   \\
%&+&A_3U_0Lsech^2(y/L)cos(k_1x) \nonumber \\
% &+&A_2U_0Lsech^2(y/L)cos(k_2x-\sigma _2t) \nonumber \\
%&+&A_1U_0Lsech^2(y/L)cos(k_1x-\sigma _1t)  \nonumber 
%\label{ROWA}
%\end{eqnarray}
%\begin{eqnarray}
%&H&(x,y,t)=c_3y-U_0Ltanh(y/L)  \\
%&+&A_3U_0Lsech^2(y/L)cos(k_1x)+A_2U_0Lsech^2(y/L) \nonumber \\
%&cos&(k_2x-\sigma _2t) +A_1U_0Lsech^2(y/L)cos(k_1x-\sigma _1t)  \nonumber 
%\label{ROWA}
%\end{eqnarray}
$H(x,y,t)=c_3y-U_0Ltanh(y/L) 
+A_3U_0Lsech^2(y/L)cos(k_1x)+A_2U_0Lsech^2(y/L) 
cos(k_2x-\sigma _2t) +A_1U_0Lsech^2(y/L)cos(k_1x-\sigma _1t) 
$.
In Fig.~\ref{FoliationFieldAndCurves}a we show  simultaneously the stable and unstable foliation fields, $f_s^t(z)$ and $f_u^t(z)$, of this system, together with curves of zero-angle Fig.~\ref{FoliationFieldAndCurves}b, $\theta(z, t)=0$, where $
\theta(z, t):=\arccos\frac{\left\langle f_s^t(z), f_u^t(z) \right\rangle}{\|f_s^t(z) \| \| f_u^t(z) \|},$ found by implicit function theorem as described in \cite{MB1}, corresponding to significant shape coherence.  
%These curves by the construction are slowly changing in shape, and hence by the theory of congruence of curves.
%, and if there is enough smoothness, regularity of curvature implies that these curves have slowly changing curvature.
The main work of this paper therefore it that we show that this detail can be skipped as the FTC we introduce significantly simplifies the geometry and facilitates the computation.

\begin{figure}[h]
%(*  \includegraphics[width=0.55\textwidth]{TestVeryGood5e8OnethresholdNum0to1em2} *)
 \includegraphics[width=0.98\linewidth, height = 0.1\textheight]{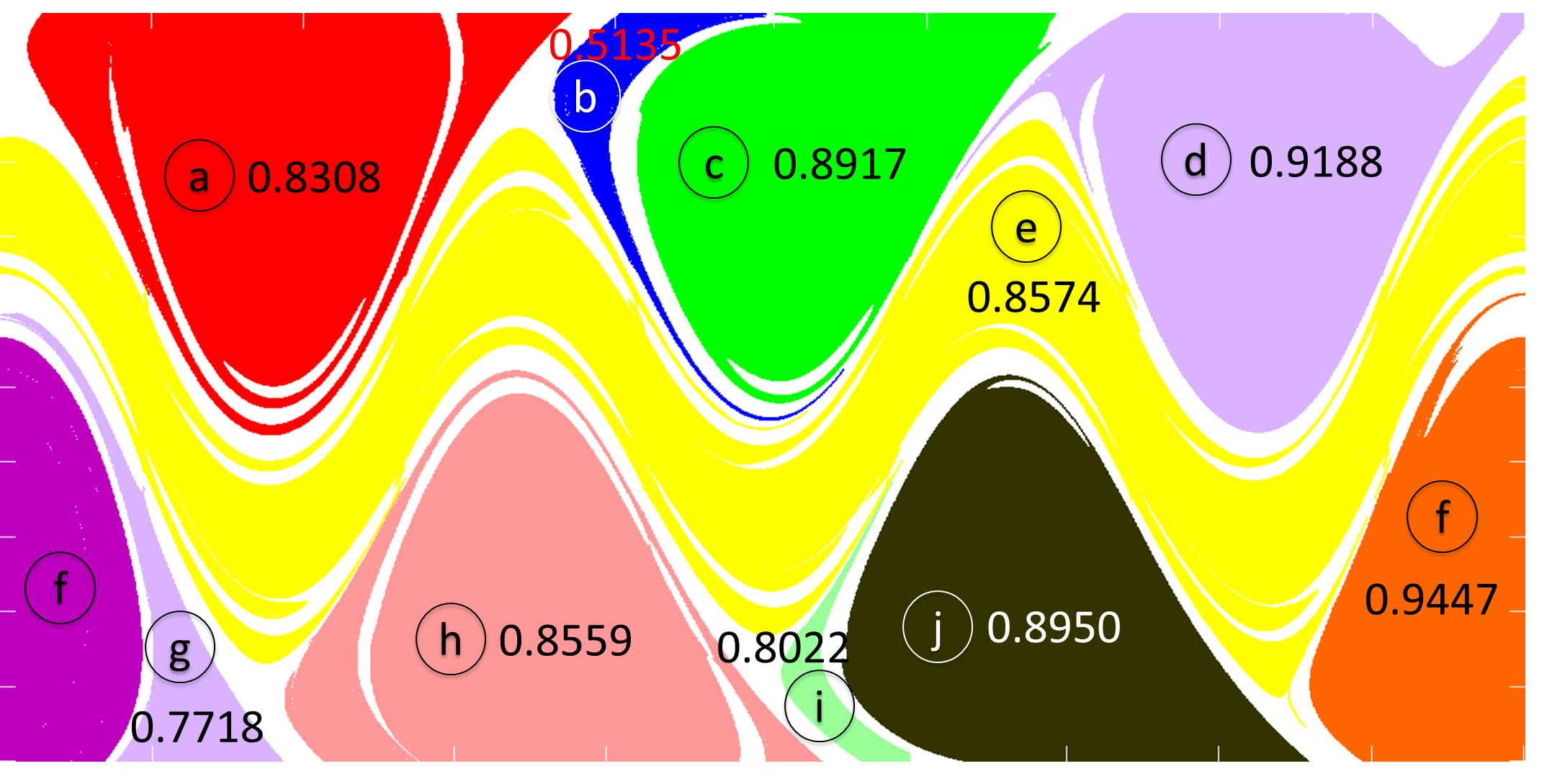} 
   \caption{The partition of the  Rossby wave system from diffusion-like ``seeded region growing" method of the FTC seen in Fig.~\ref{RWFTC}.  Numbers adjacent each letter a-j indicate the shape coherence $ \alpha(A, \Phi_T(A),T)$ of the set shown.}
 \label{Partition}
\end{figure}

\begin{figure}[h]
  \centering
  \includegraphics[width=0.48\textwidth]{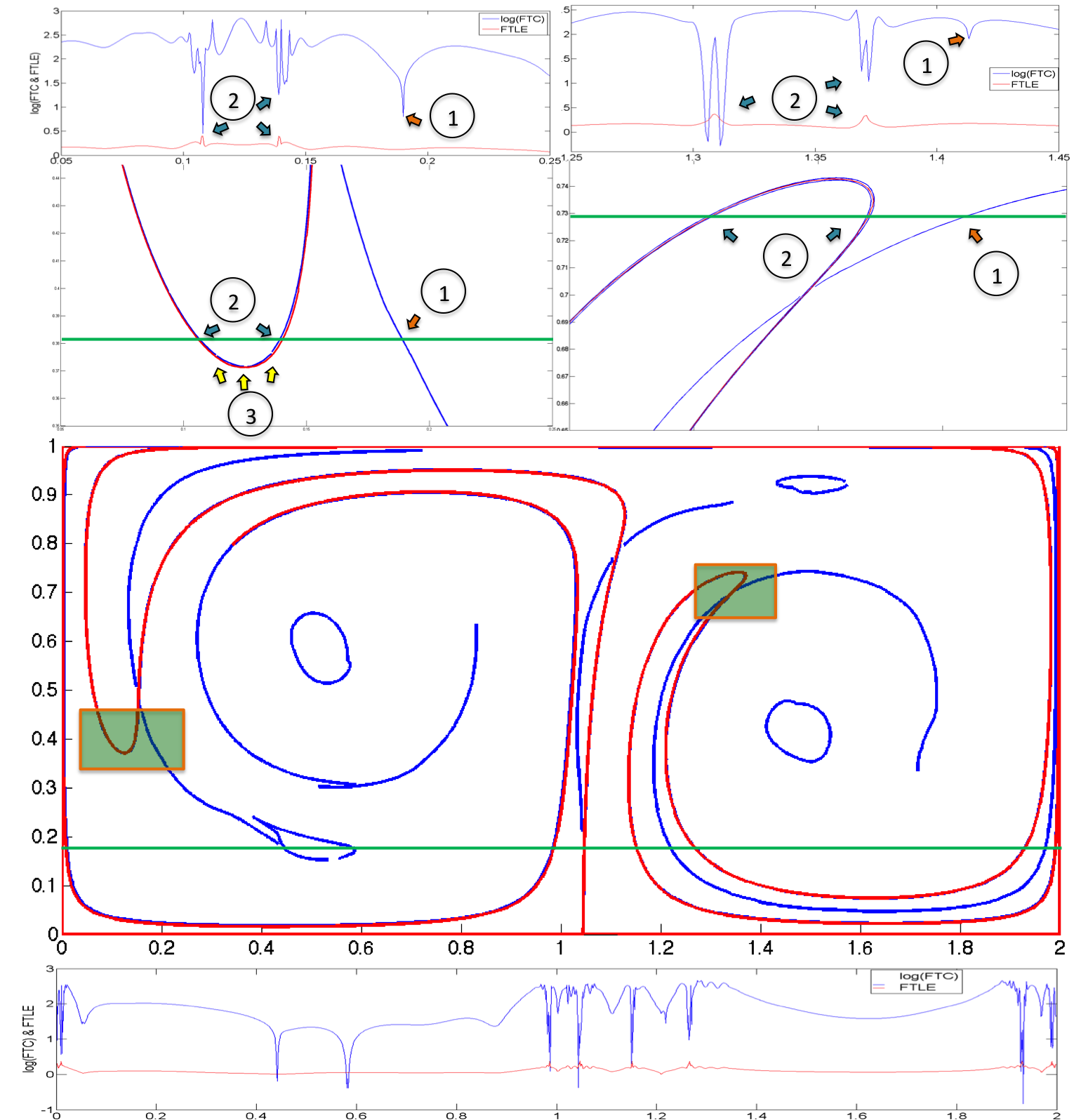} 
    \caption{ The FTC troughs of the double-gyre  highlighted in blue, and FTLE ridges highlighted in red, are seen to be sometimes in close proximity.  Thus   similar regions are indicated.  However, they are often not near each other indicating disparate regions.   In the blow-up insets above, we see labels ``2" indicate where the FTLE ridges are in-between FTC troughs, and likewise the one-dimensional slice along the green curve show these field values on a log scale repeating this outcome.  Lows of the blue curve are near but offset slightly highs of the red curve, often, highs of the blue field surrounding lows of the red field.  At such locations the two disparate computations reveal similar dynamics.  However, regions indicated by ``1" reveal FTC-trough curves which are entirely separated from any FTLE behavior of interest.  Thus a different outcome is found at such locations.  Closed curves of FTC-troughs indicate shape coherence.  At regions indicated by ``3" we see gaps in the FTC curves. Bottom we show a slice of the FTC and FTLE curve through the full phase space.}
 \label{FTC-FTLEz}
\end{figure}

%Intuition behind the FTC development is based on the 
%idea that quick detection of distorting material lines' relates to the changing curvature of those material curves passing through points where there is an especially low propensity for curvature to change. 
The intuition behind the FTC development is based the idea that  the folding behaviors  involve the maximal propensity of changing curvature.  This suggests that  regions of  space corresponding to  slowly changing curvature include  boundaries of significant shape coherence.
We define the {\bf maximum finite-time curvature (maxFTC)},  $C_{t_0}^{t_0+\tau}(z)$, and {\bf minimum finite-time curvature (minFTC)}, $c_{t_0}^{t_0+\tau}(z)$, for a point $z$  in the plane under a flow $\Phi_{t_0}^{t_0+\tau}$ over the time interval $[t_0, t_0+ \tau]$ by,
\begin{eqnarray}\label{ftc}
C_{t_0}^{t_0+\tau}(z)=\lim_{\varepsilon \to 0} \sup_{\|v\|=1} \kappa(\Phi_{t_0}^{t_0+\tau}(l_{\varepsilon, v}(z))),\\
c_{t_0}^{t_0+\tau}(z)=\lim_{\varepsilon \to 0} \inf_{\|v\|=1} \kappa(\Phi_{t_0}^{t_0+\tau}(l_{\varepsilon, v}(z)))
\end{eqnarray}
where,
$l_{\varepsilon, v}(z):=\{\hat{z}=z+\varepsilon s v, |s|<1\},
$ and $v$ is a unit vector.  So, $l_{\varepsilon, v}(z)$ is a small line segment passing through the point $z=(x,y)$, when $\epsilon<<1$. 
 Then finally we find it most useful to define the ratio of these, 
 \begin{equation}\label{ftcfield}
 r_{t_0}^{t_0+\tau}(z)=\frac{C_{t_0}^{t_0+\tau}(z)}{c_{t_0}^{t_0+\tau}(z)},
 \end{equation}
 which we simply call the {\bf finite-time curvature field,} or {\bf FTC}.
% **
%In other words, the maxFTC field is a map of  the maximal propensity of  small line segments of material curves to develop curvature  over a finite time in the flow.  
Generally when the maxFTC has a trough (curve) of small values, then this suggests that there is a 
strong nonhyperbolicity such as an elliptic island boundary or some other form of tangency as displayed in Fig.~ \ref{IdeaOfFTC}.  These are the darker blue ``FTC trough curves" we see in Fig.~\ref{RWFTC}a, and they serve as boundaries between shape coherent sets.
On the other hand, 
 the largest ridges of the maxFTC  illustrate  points where there is both significant curvature growth along one direction but small curvature growth along a transverse direction, recalling the area preservation assumption.  These level curves arise in the scenario of the sharply changing curvature developing  at the most extreme points in a (hetero)homoclinic tangle, such as illustrated in Fig.~\ref{IdeaOfFTC}.  These curves  can maintain their shape for some time.  Notice that the FTC also shows  troughs similar to the the maxFTC, but emphasized, and so these (blue) trough curves can also be used to determine shape coherent sets. A particularly interesting feature of these FTC fields is the large variation in certain regions, indicated at the top and side of Figs.~\ref{RWFTC}a; this is clearly due to co-located hyperbolicity and nonhyperbolicity regions of (hetero)homoclinic tangles, discussed in greater detail in comparison to FTLE in Fig.~\ref{FTC-FTLEz}.    The ratio FTC field  in Fig.~\ref{RWFTC}c, most clearly delineates boundaries of the shape coherent sets as low (blue) troughs.%, and  we use these to find boundaries.

To construct shape coherent sets from the FTC, we describe two complementary perspectives.  One again follows the idea of curve continuation by the implicit function theorem, but on the FTC to track a level curves of $ r_{t_0}^{t_0+\tau}(z)$.  That is, if a point $z_0$ where a (near) minimal value $ r_{t_0}^{t_0+\tau}(z)=R$ is found, representing a point in the trough, then other values nearby can be derived by $z'=h(z)=-\frac{\partial(r_{t_0}^{t_0+\tau})/\partial y}{\partial(r_{t_0}^{t_0+\tau})/\partial x}(z)  $, as an ordinary differential equation with initial condition $z(0)=z_0$, and the derivative $'=\frac{d}{ds}$ represents variation along the $s$-parameterized arc.  Furthermore, by the above regarding principle component analysis, directions of maximal curvature are also encoded in the principle vectors of $D\Phi_t(z)$.

A direct search for the interiors of sets between low troughs of the FTC  is a problem of defining regions between boundary curves and this relates to a common problem of image processing called image segmentation, \cite{ImageSegmentation}.  
%{\color{blue}There are many algorithms designed for this general problem, and we defer to the class of algorithms that depend on a seed within a region, that is if a point is identified, then  a simple diffusion process such as the head equation explores a growing set within the seeded region.  A particularly simple but robust approach is called ``seeded region growing", \cite{Adams-Bischof}, that iteratively adds neighbor pixels, checking for a threshold.  A well regarding implementation can be found at, \cite{mathworks}. }
%
%The theoretical developments describe that the slowly changing curves indicated  by the maxFTC and FTC are boundaries of shape coherent sets,  an algorithm to identify the interiors corresponding to the (blue) curves in Fig.~\ref{RWFTC}a-c must be developed.  In fact, this problem of defining regions between boundary curves as defined in an image is a now classic problem in image processing called image segmentation, \cite{**somebook on image processing}. 
%
 In particular, we applied  the diffusion-like ``seeded region growing" method \cite{Adams-Bischof} that begins with selecting a set of seed points.
Here we apply 100 uniform grid points as seeds and use 4 connected neighborhood to grow from the seed points.
We slightly improve a well regarded implementation that can be found at, \cite{mathworks}.
See Fig.~\ref{Partition} for the partitioning results. 
Several shape coherent sets corresponding to Fig.~\ref{RWFTC} are found.   Specifically, the middle yellow band has an $ \alpha(A, \Phi_T(A),T)=0.8574$, and likewise the $\alpha$-values of the rest of the colored sets are shown in Fig.~\ref{Partition}.  Some of the difference of otherwise symmetric regions are due to the region clipping as shown. 

Finally we contrast  results from the FTC field   versus the highly popular FTLE field \cite{SLM} since at first glance the pictures may seem essentially  similar, despite the significantly different definitions and different perspectives. Recall \cite{SLM} that the FTLE is defined pointwise over a time epoch-$t$ that $L_t(z)=\frac{1}{t} \log \sqrt{\rho(D\Phi_t'(z)D\Phi_t(z))}$, where $\rho$ is the largest eigenvalue of the Cauchy-Green strain tensor.  
In the following we contrast FTLE and FTC in the context of the that follows the nonautonomous Hamiltonian, $H(x,y,t)=A \cos (\pi f(x,t)) \cos (\pi y)$, where $f(x,t)= \epsilon \sin(\omega t) x^2 + (1-2 \epsilon \sin (\omega t)) x$, $\epsilon = 0.1$, $\omega=2 \pi/10$ and $A=0.1$, which has become a benchmark problem, \cite{SLM}.
Observe in Fig.~\ref{FTC-FTLEz} that sometimes an FTC trough indicative of shape coherence may occur spatially in close proximity to an FTLE ridge indicative of high finite time hyperbolicity \cite{HB, SLM} and thus suggests a transport pseudo-barrier \cite{HB, SLM}.  It is true that folding often occurs in close proximity to regions of strong hyperbolic stretching, \cite{Thiffeault,Thiffeault1, TA}, as already hinted by the fast variations of FTC in hyperbolic regions as seen in the traces on the tops and sides  of Fig.~\ref{RWFTC}ab.  However, in Fig.~\ref{FTC-FTLEz} we directly address the coincidences and differences, by locating the  troughs of the FTC shown as blue curves, and the ridges of the FTLE shown as red curves.  Clearly sometimes FTC troughs sometimes finds curves close to FTLE ridges, but sometimes entirely new curves are found.  When the FTC troughs are closed, shape coherent sets are indicated, and not found any other way when not near the FTLE.  Finally note that indicated by ``3" in  Fig.~\ref{FTC-FTLEz}, where the FTLE curves may have a strong curvature in them those FTC may be in close parallel, but the FTC trough curves may have breaks indicated.  

With these coincidences, and given differences in definitions, concepts and results, we have offered here the FTC as a new concept for interpreting shape coherence in turbulent systems, that results in a decomposition of chaotic systems into regions of simplicity, and by complement regions of complexity.  There is the promised implications that we plan to study further, between shape coherence and persistence of energy and enstrophy along Lagrangian trajectories as was likewise previously studied in the context of FTLE \cite{OG4}.

\bibliography{PaperForPRL080414}% Produces the bibliography via BibTeX.
%\printbibliography 
\end{document}